\documentclass[reqno,draft]{amsart}
\usepackage{multicol, color}
\newcommand{\rem}[1]{}

\theoremstyle{plain}
\newtheorem{lemma}{Lemma}
\newtheorem{theorem}[lemma]{Theorem}

\newtheorem{definition}[lemma]{Definition}
\theoremstyle{remark}
\newtheorem{remark}{Remark}

\newcommand*  {\R} {{\mathbb R}}
\newcommand*  {\Z} {{\mathbb Z}}

\def\Om{\Omega}

\def\pp{\partial}
\def\del{\partial}

\begin{document}

\title[Loss of smoothness and energy conservation in
the $3d$ Euler equations]
{Loss of smoothness and energy conserving rough weak solutions for
the $3d$ Euler equations}

\date{October 14, 2009}
\thanks{{\bf To appear in:} Discrete and Continuous Dynamical
Systems.}

\author[C. Bardos]{Claude Bardos}
\address[C. Bardos]
{Laboratory J.~L.~Lions\\
Universit\'{e} Pierre et Marie Curie\\ Paris, 75013, France \\
{\bf ALSO}  \\
Wolfgang Pauli Institute, Vienna, Austria}
\email{claude.bardos@gmail.com}

\author[E.S. Titi]{Edriss S. Titi}
\address[E.S. Titi]
{Department of Mathematics \\
and  Department of Mechanical and  Aerospace Engineering \\
University of California \\
Irvine, CA  92697-3875, USA \\
{\bf ALSO}  \\
Department of Computer Science and Applied Mathematics \\
Weizmann Institute of Science  \\
Rehovot 76100, Israel} \email{etiti@math.uci.edu and
edriss.titi@weizmann.ac.il}

\begin{abstract}
A basic example of  shear flow was introduced  by DiPerna and Majda
to study the weak limit of oscillatory solutions of the Euler
equations of incompressible ideal fluids. In particular, they proved
by means of this example that weak limit of solutions of Euler
equations may, in some cases, fail to be a solution of Euler
equations. We use this shear flow example to provide non-generic,
yet nontrivial, examples concerning the loss of smoothness of
solutions of the three-dimensional Euler equations, for initial data
that do not belong to $C^{1,\alpha}$. Moreover, we show by means of
this shear flow example the existence of weak solutions for the
three-dimensional Euler equations with vorticity that is  having a
nontrivial density concentrated on non-smooth surface. This is very
different from what has been proven for the two-dimensional
Kelvin-Helmholtz problem where a
 minimal regularity implies the  real analyticity of the interface.
 Eventually, we use this shear flow to
provide explicit examples of non-regular solutions of the
three-dimensional Euler equations that conserve the energy, an issue
which is related to the Onsager conjecture.

\end{abstract}
\maketitle

\vskip0.25in

{\it This paper is dedicated to Professor V. Solonnikov,
on the occasion of his  75th birthday, as token of friendship
and admiration for his contributions to research in partial
differential equations and fluid mechanics.}

\vskip0.25in

{\bf MSC Classification}: 76F02, 76B03.
\\

{\bf Keywords}: Loss of smoothness for the three-dimensional Euler
equations,  Onsager's conjecture and conservation of energy for
Euler equations, vortex sheet, Kelvin-Helmholtz.

\section{Introduction}   \label{S-1} More than 250 years after
the Euler equations have been written our knowledge of their
mathematical structure and their relevance to describe the
complicated phenomenon of turbulence is still very incomplete, to
say the least. Both in two and three dimensions certain challenging
problems concerning the Euler equations remain open. In particular,
we still have no idea of whether three-dimensional solutions of the
Euler equations, which start with smooth initial data, remain smooth
all the time or whether they may become singular in finite time. In
the case of finite time singularity it would be tempting to rely on
weak solution formulation. However, there is almost no construction,
so far,  of weak solutions for a given initial value of the
three-dimensional Euler equations. Moreover, defining an optimal
functional space in which the three-dimensional problem is {\it
well-posed  in the sense of Hadamard} is also an important issue.

Configuration where the vorticity is concentrated, as a measure, on
a curve (in $2d$) or on a surface (in $3d$) are called
Kelvin-Helmholtz flows. They seem to play a r\^ole in numerical
simulations and in the description of turbulence. However,
mathematical analysis and experiments show that in $2d$ these
configurations are extremely unstable. The main reason for this
instability being that the density of vorticity generates a
nonlinear elliptic problem (see, e.g., \cite{Duchon-Robert},
\cite{Lebeau}, \cite{wu} and references therein).

Let us observe that the conservation of energy in the $3d$ Euler
equations is always formally true. However, physical intuition and
scaling argument, i.e.~the {\it Kolmogorov Obukhov law}, lead to the
idea that non conservation of energy in the three-dimensional Euler
equations  would be intimately related to the loss of regularity.
Therefore, Onsager \cite{Onsager} conjectured the existence of a
threshold in the regularity of the   $3d$ Euler equations  that
would distinguish between solutions which conserve energy and
solutions which might dissipate energy.

For the above reasons we believe that the detailed study of explicit
examples remains extremely insightful and useful. Therefore,  this
contribution is devoted to new information that can be obtained from
the study of  the example of shear flow that was introduced by
DiPerna and Majda \cite{DipernaMajda}.

For simplicity we will consider solutions of Euler equations defined
in a domain $\Omega$ which will denote either  the whole space
$\R^3$, or the torus $(\R/ \Z)^3$  when in the latter case the
solutions are subject to periodic boundary conditions of period $1$.

Observe that when the functions   $u_1 $ and $u_3$ are smooth the
vector field
\begin{equation}
u(x,t)=(u_1(x_2), 0, u_3(x_1-t u_1(x_2))) \label{shear}
\end{equation}
is an obvious solution of the $3d$ incompressible Euler equations of
inviscid (ideal) fluids:
\begin{eqnarray} \pp_t u +
\nabla \cdot (u\otimes u)=-\nabla p\,\,\hbox{ and } \,\, \nabla\cdot
u=0\,, \label{euler}
\end{eqnarray}
with $p=0$, i.e.~this is a pressureless flow. When defined on the
torus $(\R/ \Z)^3\,$ such solutions have finite time-independent
energy, that is

\begin{equation}
\frac{d}{dt}\int_{(\R/\Z)^3}|u(x,t)|^2dx= 0\, .
\end{equation}

It is worth stressing that the following observation will be
essential for the remainder of this paper. Specifically, we observe
that the above properties remain true under much weaker assumption
on the vector field $u(x,t)=(u_1(x_2), 0, u_3(x_1-t u_1(x_2))),$
provided the notion of weak solution is used.

\begin{definition}\label{definition}  A vector field
$u\in L_{\rm loc}^2(\Omega\times [0,\infty))$ is a weak solution of
the Euler equations~(\ref{euler}) with initial data
$$
u_0\in L_{\rm loc}^2(\Omega)\,,\quad \nabla\cdot u_0=0\,,
$$
if $u$ is divergence free, in the sense of distributions in $\Omega
\times [0,\infty)$, and if for any divergence free vector field of
test functions $\phi\in C^\infty_c(\Omega \times [0,\infty)) $ one
has:
\begin{equation}
\int_{\Omega\times [0,\infty)}[u\cdot \del_t \phi +\langle u\otimes
u, \nabla \phi\rangle] dxdt= \int_\Omega u_0(x)\cdot \phi(x,0) dx.
\label{weak}
\end{equation}
\end{definition}

\begin{theorem}\label{energyconth}
(i) Let  $u_1,u_3 \in L^2_{\rm loc}(\R)$, then the shear flow
defined by $(\ref{shear})$  is a weak solution of the Euler
equations, in the sense of Definition~\ref{definition}, in
$\Omega=\R^3$.

(ii) Let  $u_1,u_3 \in L^2(\R/\Z) $ then the shear flow defined by
$(\ref{shear})$  is a weak solution of the Euler equations, in the
sense of Definition~\ref{definition}, in $\Omega=  (\R/ \Z)^3$.
Furthermore, in this case  the energy of this solution is constant.
\end{theorem}
The proof of the above statements follows from a  lemma, which is
deduced from the Fubini theorem. Below we state, without a proof,
the periodic case version of such a Lemma.

\begin{lemma}
Let $\Omega=((\R/ \Z))^3\,,\,\, u_1, u_3\in L^2((\R/ \Z))$, then for
every   test functions $\phi_i \in C^\infty(\R/ \Z)$, for $i=1,2,3$,
and $\phi_4\in C_c^\infty ([0,\infty))$ the following standard
formula
\begin{eqnarray}
&&\int_{\Om \times [0,\infty)} u_3(x_1-t u_1(x_2))\phi_1(x_1)\phi_2(x_2)\phi_3(x_3)\phi_4(t)dx_1dx_2dx_3dt\nonumber\\
&&= \int_{\Om\times [0,\infty)}u_3(x_1)\phi_1(x_1+t
u_1(x_2))\phi_2(x_2)\phi_3(x_3)\phi_4(t)dx_1dx_2dx_3dt\label{lem1}
\end{eqnarray}
is valid.
\end{lemma}

DiPerna and Majda introduced the shear flow~(\ref{shear}) in their
seminal paper \cite{DipernaMajda}  to construct a family of
oscillatory solutions of the $3d$ Euler equations whose weak limit
does not satisfy the Euler equations. In this paper we will
investigate other properties of this shear flow in order to address
issues related to the questions of well-posedness, stability of
solutions whose vorticity contains density functions that are
concentrated on surfaces (this problem being closely related to the
Kelvin-Helmholtz problem), and conservation of energy (Onsager
conjecture \cite{Onsager}). It is worth mentioning that this shear
flow was also investigated by Yudovich \cite{Yudovich} to show that
the vorticity  grows to infinity, as $t \to \infty$, which he calls
gradual loss of smoothness. This is a completely different notion of
loss of smoothness than the one presented in Theorem
\ref{loss-of-smoothness} below, where we  show the instantaneous
loss of smoothness of the solutions for certain class of initial
data.

\section{Instability of Cauchy problem and loss of
smoothness}\label{loss-smoothness-section}  Most of  the basic
existing results for the initial value problem concerning the Euler
equations~(\ref{euler}) rely on the expression of this solution in
term of the vorticity,  $\omega=\nabla \wedge u$, which satisfies in
$\R^n$, for $n=2,3$, the equivalent system (under the appropriate
boundary conditions at infinity) of equations:
\begin{eqnarray}
&&\del_t \omega + u\cdot \nabla \omega =\omega \cdot \nabla u \label{ricatti}\,,\\
&&\nabla \cdot u=0\,, \nabla \wedge u= \omega\,.  \label{biotsavart}
\end{eqnarray}
Equation (\ref{biotsavart}) defines $u$ in term of $\omega$, which
is given (in $\R^n$, for $n=2,3$) by  the Biot-Savart law; that is
$u=K(\omega)$ where $K$ is a pseudo-differential operator of order
$-1$\,. Therefore, in this case the map $\omega\mapsto \nabla u$ is
an operator of order $0$. As it is well known,  equation
(\ref{ricatti}) seems to share some similarity with the Riccati
equation
\begin{equation}
y'= Cy^2 \hbox{ whose solution is } y(t)=\frac{y(0)}{1-Ct y(0)}\, ,
\label{ricatti2}
\end{equation}
which blows up in finite time for every $y(0) > 0$. There is not
enough justification for this similarity to deduce from
(\ref{ricatti2}) some blow up property for the Euler equations.
However, one can deduce some {\it local in time} existence and
stability results in any appropriate norm $|\!|.|\!|$ which
satisfies the relation:
\begin{equation}
|\!| \omega \cdot \nabla u|\!|=|\!| \omega \cdot \nabla (K(\omega))|\!|\le C |\!|\omega|\!|^2\,.
\end{equation}
On the one hand, the operator $K$ is not continuous from $C^0$ to
$C^1$, therefore the $L^\infty$ norm is not appropriate for this
scenario. On the other hand, the H\"older norms, i.e. $\omega \in
C^{0,\alpha}$ or $u\in C^{1,\alpha}$, for $\alpha \in (0,1]$, are
convenient. With the standard Sobolev estimates the norm $H^s$, for
$s>\frac52$, i.e. $\omega \in H^{s-1}$ or $u\in H^s$, would  also be
convenient (and leads, by virtue of common functional analysis
tools, to slightly simpler proofs, see, e.g.~\cite{Majda-Bertozzi}).
This is fully consistent with the fact that $H^s$, for $s>\frac52$,
is continuously imbedded in $C^{1,s-\frac52}$.

With this classical observations in mind we recall the following
facts (see also the recent surveys for more details
\cite{Bardos-Titi} and \cite{Constantin}):
\begin{enumerate}

\item  For initial data $u(x,0)=u_0(x)$ in $C^{1,\alpha}$ the Euler equations
(\ref{euler}) has a unique local in time solution $u(x,t)$ in
$C^{1,\alpha}$ (cf. \cite{lichtenstein}). The same result is valid
for initial data in $H^s$, for $s>\frac 52$ (cf.
\cite{BealeKatoMajda},\cite{Majda-Bertozzi}). Moreover, this unique
solution conserves the energy.  In spite of the fact that the above
results imply the short time control of the  $L^\infty$ norm of the
vorticity (which seems to be the relevant quantity) one has the
following complementary statement established in
\cite{BealeKatoMajda} (see also \cite{Majda-Bertozzi}). For every
initial data $u(x,0)$ in $C^{1,\alpha}$ or in $H^s$, for
$s>\frac52$, the solution of the three-dimensional Euler equations
exists and depends continuously on the initial data, for  as long as
the time integral of the $L^\infty$ norm of the vorticity remains
bounded.

\item  Following \cite{DeLellisandSzekelyhidi} one
can prove  (in any space dimension) the existence of initial data
$u_0\in L^2(\Om)$ (not explicitly constructed) for which the Cauchy
problem has, with the same initial data, an infinite family of weak
solutions of the  Euler equations: a residual set in the space $C
(\R_t; L_{\rm weak}^2(\Om))\,.$

\item Eventually one does not know the existence of a $3d$ regular
(say in $C^{1,\alpha}$) solution of the Euler equations that becomes
singular in a finite time (blow up problem).

\end{enumerate}

The shear flow (\ref{shear}) has also been  used by DiPerna and
Lions  \cite{Lions}, and  by \cite{Yudovich}, as an example to
demonstrate some issues related to the instability of the solutions
of the three-dimensional Euler equations. In particular,  DiPerna
and Lions (cf. \cite{Lions} page 124) have established the following

\begin{theorem} [DiPerna-Lions]\label{dipernalions}
For every $p \ge 1$, $T >0$ and $M > 0$ given there exists a smooth
shear flow solution of the form (\ref{shear}) for which
$\|u(x,0)\|_{W^{1,p}}=1$ and $\|u(x,T)\|_{W^{1,p}}> M$.

\end{theorem}
In fact the proof of this theorem that has been  presented in
\cite{Bardos-Titi} (Proposition 3.1) shows that, for every $p \ge
1$, there exist shear flow solutions of the form (\ref{shear}) with
$ u(x,0)\in  W^{1,p} $ and $u(x,t)\notin W^{1,p}$ for any
$t\not=0\,.$ Here, we show, in addition,  the instantaneous loss of
smoothness of weak solutions for the $3d$ Euler equations with
initial data in the H\"older space $C^{0,\alpha}$, with $\alpha \in
(0,1)$. This underlines the r\^ole of the space $C^1$ as the {\it
critical space} for short time well-posedness of the $3d$ Euler
equations; namely: for initial data {\it more regular than $C^1$},
say in $C^{1,\beta}$, with $\beta \in (0,1]$,  one has
well-posedness  for the $3d$ Euler equations, and for {\it less
regular} initial data, specifically initial data in $C^{0,\alpha}$,
with $\alpha \in (0,1)$, one has ill-posedness.

\begin{theorem}\label{loss-of-smoothness}
(i) For $u_1(x),u_3(x)\in C^{1,\alpha}$, with $\alpha \in (0,1]$,
the shear flow solution (\ref{shear}) is in $C^{1,\alpha}$, for all
$t\in \R$.

(ii) For $u_1(x),u_3(x)\in C^{0,\alpha}$, with $\alpha \in (0,1)$,
the shear flow solution (\ref{shear}) is always in $C^{0,\alpha^2}$.

(iii) There exist shear flow solutions, of the form (\ref{shear}),
which for $t=0$ belong to $C^{0,\alpha}$, for some $\alpha \in
(0,1)$, and which for $t\not=0$ do not belong to $C^{0,\beta}$ for
any $\beta>\alpha^2$.
\end{theorem}

{\bf Proof}. Observe first that in (i) the evolution of regularity
concerns only the component $u_3$ ($u_1$ remains $t$ independent).
The statement (i) is trivial, but it is worth noticing as it shows
that our analysis is in line with the classical results of
\cite{lichtenstein}. To prove  (ii)  we write
\begin{eqnarray}
&&\frac{|u_3(x_1-tu_1(x_2+h))-u_3(x_1-tu_1(x_2))|}{h^{\alpha^2}}\nonumber\\
&&= \frac{|u_3(x_1-tu_1(x_2+h))-u_3(x_1-tu_1(x_2))|}{|tu_1(x_2+h) -tu_1(x_2)|^\alpha}\Bigg(\frac {|tu_1(x_2+h) -tu_1(x_2)|}{ h^\alpha}\Bigg)^\alpha\nonumber\\
&&\le|t|^\alpha
|\!|u_3|\!|_{0,\alpha}|\!|u_1|\!|_{0,\alpha}^\alpha\,.
\end{eqnarray}

For the point (iii) of the statement one introduces two periodic
functions $u_1(\xi)$ and $u_3(\xi)$ which near the point $\xi=0$
coincide with the function $|\xi|^\alpha$. Consequently, for every
given $t$  and  for  $x_1$ and $x_2$ small enough,
$u_3(x_1-tu_1(x_2))$ coincides with the function
$$|x_1-t|x_2|^\alpha|^\alpha .$$
In particular, for $t$ given, and for  $(x_1,x_2,x_3)=(0,x_2,x_3)$,
with $x_2$ small enough,  one has
$$u_3(x_1-tu_1(x_2))= |t|^\alpha|x_2|^{\alpha^2},$$
and the conclusion follows.

\begin{remark}\label{Lemarie}
P.G.~Lemari\'e-Rieusset observed (private communication)
 that the criticality  aspect of the space $C^1$, in the above
context,  can be sharpened by considering the Besov and the
Triebel-Lizorkin spaces, $B^s_{p,q}$ and $F^s_{p,q}$, respectively.
Indeed, one has, on the one hand,  the inclusions (see, e.g.,
\cite{Lemarie-Rieusset} and \cite{Triebel})
\begin{equation}
C^{1,\alpha}= B^{1+\alpha}_{\infty,\infty}\subset
B^{1}_{\infty,1}\subset C^1\subset  F^1_{\infty,2} \subset
B^{1}_{\infty,\infty}\subset B^{
\beta}_{\infty,\infty}=C^{0,\beta}\,,
\end{equation}
for all $\alpha\in (0,1]$ and $\beta \in (0,1)$. On the one hand,
the short time  well-posedness of the $3d$ Euler equations has been
recently proven in the space $B^{1}_{\infty,1}$ by Pak and Park
\cite{Pakpark}, and on the other hand, calculations inspired by the
above proof lead to the construction of shear flows $u(x,t)$, of the
form~(\ref{shear}), which satisfy:
\begin{eqnarray*}
u(x,0)\in   F^1_{\infty,2} \hbox{ and for all } t\not=0
\,\,u(x,t)\notin F^1_{\infty,2}; \hbox{ or } u(x,0)\in
B^{1}_{\infty,\infty} \hbox{ and for all } t\not=0 \,\,u(x,t)\notin
B^{1}_{\infty,\infty}\,.
\end{eqnarray*}
The details of this analysis and further applications will be
reported in a forthcoming paper.
\end{remark}
\section{ Vorticity Surface density for shear flow and the Kelvin-Helmholtz problem}\label{khsf}
  For  $2d$ Euler equations existence of a weak solutions, with
  a given single signed Radon measure initial data
  for the vorticity density which is supported on a
curve,  has been established by Delort \cite{DE}.
  The condition on the sign of the vorticity   has
  later  been slightly relaxed \cite{LNZ}. There is no
  such theorem  in $3d$ case and the only available
  result  for initial data having a density of
  vorticity {\it concentrated } on a surface is a local in
  time  existence and uniqueness result  under very restrictive
  analyticity hypothesis
  (see \cite{Sulemsulembardosfrisch}). As a consequence
  it may be interesting to exhibit two examples of shear
  flows with nontrivial surface density that would emphasize
  the difference between the $2d$ and $3d$ situations.

{\bf Example 1}

To present  a vorticity concentrated on a surface for the shear
flow~(\ref{shear}) then it has to be of the  following form
 \begin{equation}
u_1(s)=\left\{\begin{array}{cl} & \alpha_1\,\, \hbox{ for } \,\,s<\xi_2\,\\
& \beta_1\,\, \hbox{ for } \,\, s>\xi_2
\end{array}
\right. \hskip 0.1in \hbox{ and } \hskip 0.1in
u_3(s)=\left\{\begin{array}{cl} & \alpha_3\,\, \hbox{ for } \,\, s<\xi_1\,\\
& \beta_3\,\, \hbox{ for } \,\, s>\xi_1
\end{array}
\right. ,
\\  \nonumber
\end{equation}
for some fixed real parameters
$\alpha_1,\alpha_3,\beta_1,\beta_3,\xi_1,\xi_2$, satisfying
$\alpha_1 \ge \beta_1$ and $\alpha_3 \neq \beta_3$.
The vorticity is therefore   concentrated on
the singular surface:
\begin{eqnarray*}
 \Sigma(t)=\{(x_1,x_2,x_3)| \, x_2=\xi_2\}\cup
  \{(x_1,x_2,x_3)| \,
x_1= \xi_1+ t \alpha_1, \, x_2\le\xi_2\}
 \cup \{(x_1,x_2,x_3)| \,
x_1= \xi_1+ t \beta_1, \, x_2\ge\xi_2\}\,.
\end{eqnarray*}

{\bf Example 2 }

In $3d$ with   the following configuration
 \begin{equation}
u_3(s)=\left\{\begin{array}{cl} &1\,\, \hbox{ for } x<0\,\\
& 0\,\, \hbox{ for } x>0
\end{array}
\right. \,,
 \end{equation}
 and $y=u_2(s)$ is a $C^1$ curve, the shear flow:
 $$
 u(x)= (u_1(x_2), 0, u_3(x_1-tu_1(x_2)))
 $$
 is a weak solution of the $3d$ Euler equations with a singular
 vorticity which is concentrated  on the surface
 $$
 \Gamma(t)= \{ (x_1, x_2, x_3) | ~~x_1=t u_1(x_2)\}
 $$
and is given by:
 \begin{equation}
 \omega(x,t)= (-t \frac{\del_{x_2} u_1}{(|t\del_{x_2} u_1|^2+1)^\frac12}\otimes \delta_{\Gamma(t)},
 \frac{1}{(|t\del_{x_2} u_1|^2+1)^\frac12}\otimes \delta_{\Gamma(t)}, -\del_{x_2} u_1(x_2))\,. \nonumber
 \end{equation}

The discussion  below, concerning the difference between the $2d$
and $3d$ Kelvin-Helmholtz problem, is motivated, among other things,
by the following remark.

\begin{remark} \label{energy} Example 1 is a solution of
the $3d$ Euler equations with a density of vorticity concentrated on
a surface with corners. It is  unknown whether the construction of
the same type of configuration is possible in $2d$ case. In Example
2 the function $x_2\mapsto u_1(x_2)$  does not need to be more
regular than $C^1$ to sustain a $C^1$ vorticity surface density of
the corresponding shear flow solution of the $3d$ Euler equations.
Moreover, no matter how regular this surface is initially  its
regularity will be preserved by the dynamics. Furthermore, and by
virtue of Theorem \ref{energyconth}, both examples are weak
solutions of the $3d$ Euler equations, and when considered in torus
$(\R/\Z)^3)$ they both conserve energy.
 \end{remark}

In an attempt to  understand  the effect of the dimension  it seems
appropriate to compare Example 2 with classical results concerning
the Kelvin-Helmholtz problem.


As we have mentioned in the introduction the Kelvin-Helmholtz
problem corresponds to the situation where the vorticity is
concentrated on a moving  orientable curve in $r(t, \lambda)$, in
$2d$, parameterized by a parameter $\lambda \in \R$, or on a moving
orientable surface $r(t,\lambda, \mu)$, in $3d$, parameterized by
the parameters $(\lambda,\mu)\in \R^2$.

We assume  that the curves or the surfaces are $C^1$  orientable
manifolds, denoted by $\Gamma(t)$, with unit normal $\vec n$. For
$x\notin \Gamma(t)$, the velocity $u$ can be expressed explicitly in
term of the vorticity by the following Biot-Savart formulas:
\begin{equation}
 u(x,t)=\left\{\begin{array}{cl} \frac1{2\pi}R_{\frac{\pi}{2}}  \int \frac {x-r(t,\lambda')}
{|x-r(t,\lambda')|^2} \tilde \omega(t,r(t,\lambda')) |\del_\lambda r(t,\lambda')|d\lambda' ~~~~~\hbox{ in } 2d  \,,  \vspace{4mm}\\
 -\frac1{4\pi}   \int \frac {x-r(t,\lambda',\mu')}
{|x-r(t,\lambda',\mu')|^3} \tilde \omega(t,r(t,\lambda',\mu'))
|\del_\lambda r (t,\lambda',\mu') \wedge\del_\mu r (t,\lambda',\mu')
|d\lambda'd\mu' ~~~~~\hbox{ in }  3d\, , \end{array} \right.
 \\  \label{bios}
\end{equation}
where $R_{\frac{\pi}{2}}$ is  the $\frac{\pi}{2}$ rotation matrix,
and $\tilde \omega$ is the {\it vorticity density} on these
manifolds.

When $x$ converges to a point $r \in \Gamma(t)$ the velocity
$u(x,t)$ converges to two different values, on either side of the
manifold, $u_\pm(r,t)$. In particular, and  in agreement with the
divergence free condition, one has
\begin{equation}
u_+(r,t)\cdot \vec n=u_-(r,t)\cdot \vec n\,\,\,,\quad \omega(x,t)=
(u_+(r,t)-u_-(r,t))\wedge \vec n \otimes \delta_{\Gamma(t)}(x)\,,
\end{equation}
for $r \in \Gamma(t)$ and $x \in \R^d$, $d=2,3$.

The  vorticity density $\tilde \omega$ is a vector valued density.
In the $2d$ case this vector is orthogonal to the plane of the flow
and therefore is identified with a scalar. Hence, the vorticity
density is related to the vorticity by the expressions:
\begin{eqnarray} \omega(x,t)& = &(u_+(r,t)-u_-(r,t))\wedge \vec n
\otimes
\delta_{\Gamma(t)}(x) \nonumber \\
& = & \left\{
\begin{array}{cl}
\tilde \omega(t,r(t,\lambda)) |\del_\lambda r(t,\lambda)|d\lambda ~~~~~~\hbox{ in } 2d\,, \vspace{2mm}\\
\tilde \omega(t,r(t,\lambda,\mu)) |\del_\lambda r
(t,\lambda,\mu)\wedge\del_\mu r (t,\lambda,\mu) |d\lambda d\mu\,
~~~~~~\hbox{ in } 3d\,. \label{BS}
\end{array}
 \right.
 \end{eqnarray}
Formulas (\ref{bios})  remain valid for $x\in \Gamma(t)$ with the
integral taken in the sense of Cauchy principal value and with the
left-hand side of (\ref{bios}) replaced by the averaged velocity
\begin{equation}\label{bios2}
v=\frac{u_++u_-}2\,.
\end{equation}
Therefore, with  some hypothesis on the regularity of the solution
(cf. \cite{Lopes} for details)  the problem can be reduced to
equation (\ref{bios}) for $v$ with:
\begin{equation}
(\del_t r-v)\cdot \vec n =0\,, \label{lagrange}
\end{equation}
and in $2d$
\begin{equation}
 \del_t \tilde \omega  +\frac{\del}{\del \lambda}\Bigg( \frac{\tilde \omega}{|r_\lambda|^2}(v-r_\lambda)\cdot r_\lambda \Bigg)=0\,\label{dynamique}
 \end{equation}
 or   in $3d$ with $N= \del_\lambda r(t,\lambda, \mu) \wedge \del_\mu  r(t,\lambda, \mu) $
\begin{eqnarray}\label{3d-dynamic}
 &&\del_t \tilde \omega  +\frac{\del}{\del \lambda}\Bigg( \frac{\tilde \omega}{|N|^2}(((v-\del_t r) \wedge \del_\mu r) \cdot N )\Bigg)
 -\frac{\del}{\del \mu}\Bigg( \frac{\tilde \omega}{|N|^2}(((v-\del_t r) \wedge \del_\lambda r) \cdot N )\Bigg )
\nonumber\\
&& = \frac{1}{|N|^2}((\del_\mu r\wedge N) \cdot \tilde \omega)
\del_\lambda v-
 \frac{1}{|N|^2}((\del_\lambda r\wedge N) \cdot \tilde \omega))\del_\mu  v  \,.
  \end{eqnarray}

 We recall below some classical results which contribute to the understanding
 of the basic properties of this
 problem (see also, e.g., \cite{Bardos-Titi}).

\begin{enumerate}

\item  The initial value problem is locally, in time, well-posed in
both, the $2d$ and the $3d$, cases in the class of analytic data.
More precisely, for any   initial   curve (respectively  surface)
$\Gamma(0,\lambda)$, (respectively $\Gamma(0,\lambda,\mu)$)  and any
initial density of vorticity $\tilde \omega(0,r(0,\lambda))$
(respectively $\tilde \omega(0,r(0,\lambda,\mu))$ which can be
extended as analytic functions uniformly bounded in the strip $ |\Im
\lambda| \le c$, in the complex plane $\lambda\in \mathbb{C}$, for
some $c>0$, (respectively $ |\Im \lambda| + |\Im \mu|\le c$, for
$(\lambda,\mu)\in \mathbb{C}^2$, and for some $c>0$) there exists a
finite time $T$ and a constant $C$ such that the initial value
problem~(\ref{lagrange}) and~(\ref{dynamique}) (respectively,
(\ref{lagrange}) and~(\ref{3d-dynamic})) has, for $ 0 \le t <T$, a
unique solution which is analytic in the strip $ |\Im \lambda| \le
C(T-t) $ (respectively $ |\Im \lambda| + |\Im \mu|\le C(T-t)$) (cf.
\cite{Sulemsulembardosfrisch}).

\item  There exist in $2d$ (to the best of our knowledge this
 issue has not been addressed in $3d$) analytic solutions that become
 singular in finite time. This has been first
observed by numerical simulations of Baker, Meiron and Orszag
\cite{Bakermeironorszag}, then Duchon and Robert
\cite{Duchon-Robert} have shown the existence   of  a very large
class of singularities which can be reached in a finite time by
analytic solutions. Eventually, Caflisch and Orellanna
\cite{caflishorellana} have constructed analytic solutions, for $0
\le t<T$, which exhibit a cusp as $t$ approaches $T$. Specifically,
with $0<\nu<1$ they have shown that their solutions satisfy:
\begin{equation}
\lim_{t\rightarrow T}(\Gamma(t,\cdot), \tilde
\omega(t,r(t,\cdot)))=(\Gamma(T,\cdot),\tilde\omega(T, r(T,\cdot)) )
\left\{\begin{array}{cl}
& \notin C^{1,\nu}\times C^\nu\,,\\
&\in C^{1,\nu'}\times C^{\nu'} \,\,\,\, \hbox{ for every} \,\,\,
\nu'\in (0,\nu)\,.
\end{array}
\right.
 \\ \nonumber
 \end{equation}
 \item In $2d$ : If in a $(t,\lambda)$ neighborhood of a point
 $(t_0,\lambda_0)$  the  vorticity density,
 $\tilde\omega(t,r(t,\lambda))$, does not vanish and if
 the functions $r(t ,\lambda ), \tilde\omega(t , r(t,\lambda))$ have some {\it limited  regularity}   then in fact they are analytic in this neighbourhood.
 By a limited regularity we mean, for instance, that in
 this neighborhood
 \begin{eqnarray}
&&(r(t ,\cdot ), \tilde\omega(t ,r(t,\cdot )) \in C^{1,\alpha}\times C^\alpha\label{lebeau} \\
&& |\lambda-\lambda'| \le C |r(t ,\lambda )-r(t ,\lambda ')|,
\,\hbox{ with some constant} ~~~ C<\infty . \label{chordarc}
\end{eqnarray}
The hypothesis (\ref{chordarc}) is called the {\it chord-arc
property},  and the  hypothesis  (\ref{lebeau}) matches perfectly
the example studied in \cite{caflishorellana} . In fact under the
chord-arc hypothesis a refined version of this statement has been
obtained by  Wu \cite{wu}, which matches some numerical observations
made by Krasny \cite{krasny}. The consequence of this observation is
that solutions {\it with limited regularity} do not exist in $2d$.
That is, if at some time $t_0$ and at some point $\lambda_0$ the
solution, $(r(t,\lambda),\tilde\omega(t ,r(t,\lambda )) $,  ceases
to be analytic then it cannot be  of limited regularity at a later
time. For instance the solution of \cite{caflishorellana} is no
longer in $C^{1,\nu'}\times C^{\nu'}$, for any $\nu'>0$, for $t>T$.
 \end{enumerate}

 \begin{remark}
 The hypothesis that $\tilde\omega(t,r(t,\lambda))$ does not vanish
 is natural. This is because if $\tilde\omega  $ vanishes near
 $(t_0,\lambda_0)$ then there is no more interface, and the
 ellipticity as described below is lost. This will appear
 explicitly in  formulas (\ref{nl1}) and (\ref{nl2}) below.
 \end{remark}


The clue in the above $2d$ results, which have been described under
different forms in \cite{Duchon-Robert}, \cite{Lebeau} and
\cite{wu}, lies in the fact that under the above hypothesis the
problem is locally a {\it small } perturbation of a linear elliptic
system. Indeed,  since this
  analysis is local one can assume, without loss of generality,
  that $\Gamma(t)= (x, \epsilon y(x,t))$ is a graph. As a result,
   equations (\ref{bios2}) , (\ref{lagrange}) and (\ref{dynamique}) are equivalent to the system:
\begin{eqnarray}
&&\del_t y- v_2=  ( v_1\del_x y)\,, \label{t3} \\
&&\del_t\tilde \omega + \del_x(v_1  \Omega_0)= -\epsilon\del_x(v_1 \tilde \omega)  \,,\label{v3}\\
&& v_1(x,t)=-\frac1{2\pi}P.V.\int \frac{y(x,t)-y(x',t)}{(x-x')^2+\epsilon^2(y(x,t)-y(x',t))^2}(\Omega_0+\epsilon \tilde \omega)dx'\,,\label{bs3}\\
&& v_2(x,t)=\frac1{2\pi}P.V.\int \frac{x-x'}{(x-x')^2+\epsilon^2(y(x,t)-y(x',t))^2}(\Omega_0+\epsilon \tilde \omega)dx'\label{bs4}\,.
\end{eqnarray}
For small values of $\epsilon$, this  system describes a small
perturbations in $\R^2$ about the stationary solution
$$
y(x,0)=0\,, u_- = \frac{\Omega_0}{2}\,,
u_+=-\frac{\Omega_0}{2}\,.
$$
Indeed, for functions $f$ and $y$ in
$C^1$, with $\frac{\partial y}{\partial x}$ bounded, the expansion
\begin{eqnarray}\label{expansion}
&& \frac1{ \pi}P.V.\int \frac{f(x)-f(x')}{(x-x')^2+\epsilon^2(y(x,t)-y(x',t))^2} dx'=\nonumber\\
 &&\frac1{ \pi}P.V.\int \frac{f(x)-f(x')}{(x-x')^2}\Bigg( 1 + \sum_{n\ge 1} (-1)^n\epsilon^{2n}\Bigg(\frac {y(x)-y(x'))}{x-x'}\Bigg)^2\Bigg) dx'
\end{eqnarray}
  leads to the introduction of the operators (Hilbert transform):
\begin{eqnarray}
&&Hf(x)=\frac{1}{\pi} P.V. \int\frac{1}{x-x'} f(x')dx'= \mathcal{F}^{-1}(-i {\mathrm {sgn} }(\xi) \hat f (\xi))\,\label{ht1}\\
&&|D|f(x)=\frac1{ \pi}P.V.\int \frac{f(x)-f(x')}{(x-x')^2} = \del_x
(Hf(x))=  \mathcal{F}^{-1}(|\xi| \hat f (\xi))\,.\label{ht2}
\end{eqnarray}
 This in turn gives, together with  formulas
 (\ref{t3})-(\ref{expansion}), for the perturbation about the
 stationary solution the system:
 \begin{eqnarray}
&&\del_ty_{x}- \Omega_0|D|\tilde{\omega}= \epsilon F(y_x,\tilde{\omega})_x\nonumber\\
&&\del_t \tilde{\omega} - |D| y_x= \epsilon G(y_x,\tilde{\omega})_x
\,,\nonumber
\end{eqnarray}
where in  the right-hand side  $F$ and $G$ are first order
operators. Eventually  with the introduction of the ``Laplacian" one
has:
\begin{eqnarray}
&&\del_{tt} (y_x) +\Omega_0^2\del_{xx} (y_x)= \epsilon(\del_t(F(y_x,\tilde{\omega})_x)+ |D|(\epsilon G(y_x,\tilde{\omega})_x) \label{nl1}\,, \label{elliptic1}\\
&&\del_{tt} (\tilde{\omega}) +\Omega_0^2\del_{xx} (\tilde{\omega})=
\epsilon(|D|(F(y_x,\tilde{\omega})_x)+ \del_t(\epsilon
G(y_x,\tilde{\omega})_x) \label{nl2}\,. \label{elliptic2}
\end{eqnarray}

We remark that Example 2 is not an exact solution of the $3d$
Kelvin-Helmholtz problem due to the fact that in this case the
function
$$\del_{x_1}u_3(x_1)$$
is not, as in the Example 1, a Dirac mass. However, we {\it
conjecture}, and that may be the object of future contribution, that
a solution of the $2d$ Euler equations with a vorticity of the form
\begin{equation}
\nabla\wedge u(.,t)= \omega_1(t) \otimes\delta_{\Gamma(t)} + \omega_2(t)
\end{equation}
with $r(t ,\lambda ),  \omega(t , r(t,\lambda))$ having  some {\it
limited  regularity} in the above sense and $\omega_2 \in
C^{1+\alpha}(\Omega\times\R_t)$ will exhibit the same type of
smoothing effect as in the case of the $2d$ Kelvin-Helmholtz. For
instance under these hypothesis the surface $\Gamma(t)$ should
belong to $C^\infty$, or even analytic. The intuition for this
conjecture stems from the fact that equation (\ref{BS}) is modified
by the addition of lower order terms, hence the conclusions are
expected to be similar. Now for the Example 2; this regularity
property is not true for the surface
$$
 \Gamma(t)= \{ (x_1, x_2, x_3) |~~ x_1=t u_1(x_2)\}\,.
 $$
The reason for the difference would be that  in $2d$ the
 smoothing effect  is due to the ellipticity  of the
 linearized operator while in $3d$ the situation is
 different as follows: As it was done in $2d$ case, we
 consider a local perturbation about the stationary solution.
 In this situation  we  assume
 (following the notation of   \cite{Sulemsulembardosfrisch} or
 \cite{Chandrashekar}) that  $\Gamma(t)$ can be parameterized in
 the form $x_3=\epsilon x(x_1,x_2,t)$, and reduce the analysis
 to the properties of the {\it small} perturbation about the
 stationary state
 $x_3=0, \tilde\omega^0 (x_1,x_2)= ( \tilde\omega_1^0, \tilde\omega_2^0,0)$.
 The leading part of the perturbed equations
 (as was done above in the $2d$ case) is the linear operator
 (written in the $2d$ Fourier variables $k=(k_1,k_2)$, the
  dual of $(x_1,x_2)$)
 \begin{equation}
 \del_t \left(\begin{array}{c} \hat {x}_3\\   \hat {  \omega}_1 \\ \hat{  \omega}_2\\ \hat{  \omega}_3\end{array}\right)={\mathcal A}\left(\begin{array}{c} \hat {x}_3\\   \hat {  \omega}_1 \\ \hat{  \omega}_2\\ \hat{  \omega}_3\end{array}\right)
 \end {equation}
 where
 \begin{equation}
 {\mathcal A}=\left(\begin{array}{clcr} 0 &\frac i2 \sin \theta  & -\frac i2 \cos \theta  & 0 \\
- \frac i2 |k|^2|\tilde\omega^0|^2\sin \theta &0  &0 &\frac12(k\cdot\tilde\omega^0)\sin \theta\\
\frac i2 |k|^2|\tilde\omega^0|^2\cos \theta &  0 &0 &-\frac12(k\cdot\tilde\omega^0)\cos \theta\\
0 &-\frac12(k\cdot\tilde\omega^0)\sin \theta
&\frac12(k\cdot\tilde\omega^0)\cos \theta  & 0
\end{array}\right)\,,
 \end{equation}
 with $k=(k_1,k_2)= |k| (\cos \theta, \sin \theta)$.

The eigenvalues of the matrix ${\mathcal A}$ are
$$
\{0,0,-\frac12|k\wedge \tilde\omega^0|, \frac12|k\wedge
\tilde\omega^0|\}\,.
$$
Therefore, the first order pseudo-differential operator
$$
\del_t -{\mathcal A}
$$
is no longer elliptic, as the situation is in the $2d$ case (see
(\ref{elliptic1})-(\ref{elliptic2})).

\section{Energy conservation for rough solutions}

It has been conjectured by Onsager \cite{Onsager} that for some weak
solutions of the $3d$ Euler equations the decay in energy
would be related to some loss of regularity in these solutions.
Arguing by some dimensional analysis, the H\"older exponent $1/3$
appears to be a critical value of such regularity.

On the one hand, it has been shown rigorously in \cite{Constantin-E-Titi}
that the formal conservation of energy in the $3d$ Euler Equations
is in fact true for any weak solution which is slightly more regular
than the Besov space ${\mathcal B}^{\frac13}_{3,\infty}$ (see also
 \cite{Cheskidov-Constantin-Friedlander} and \cite{Eyink}).
On the other hand, the existence of very weak solutions {\it wild
solutions} that become identically $0$ after a finite time has been
established in \cite{scheffer}, \cite{shnirelman}  and most recently
in \cite{DeLellisandSzekelyhidi}. Moreover, it is commonly believed
that for solutions which are slightly  weaker  than  ${\mathcal
B}^{\frac13}_{3,\infty}$ there might  be no conservation of energy.
In fact Eyink \cite{Eyink} has constructed a function $u_0(x)\in
C^{0,\frac13}$ which cannot be the initial data of any weak solution
which conserves the energy. This, however,  is not a {\it complete}
counter example because the existence of weak solutions for the $3d$
Euler equations with such initial data is still an open problem.

With the shear flow solution of the $3d$ Euler equations:
\begin{equation*}
u(x,t)=(u_1(x_2), 0, u_3(x_1-t u_1(x_2)))
\end{equation*}
in the torus $(\R/\Z)^3$,  it follows from  Theorem
\ref{energyconth}, with $u_1, u_3 \in L^2(\R/\Z)$, that there is no
hope for a general theorem stating that the conservation of energy
implies some type of regularity.

Observe that the hypothesis on the initial data here are much
weaker
than those for which the Onsager conjecture is stated in
\cite{Constantin-E-Titi},   \cite{Eyink}  or \cite{shvydkoykh} (see
also \cite{Cheskidov-Constantin-Friedlander}).

In \cite{shvydkoykh} Shvydkoy considers the energy conservation for
weak solutions of the Euler equations with singularities on a curve
(in $2d$) and on a surface (in $3d$). This class of solutions
includes the Kelvin-Helmholtz problem  discussed in  section
\ref{khsf}. In fact the results in \cite{shvydkoykh} turn out to be
more relevant for the Kelvin-Helmholtz problem in the
three-dimensional case rather than in two-dimensional one. The
reason being, as we have mentioned in section \ref{khsf} above, that
in the $2d$ case a minimal regularity for the Kelvin-Helmholtz
problem implies analyticity; and therefore the conservation of
energy of the solutions follows, while in the $3d$ case  the
ellipticity of the linearized operator is no longer true and there
is room for less regular (non-analytic), and possibly singular,
surface solution of the $3d$ Kelvin-Helmholtz problem. In agreement
with this observation we propose the following example. Consider in
Theorem \ref{energyconth} the shear flow~(\ref{shear}) in the torus
$(\R/\Z)^3$, with $u_1, u_3 \in L^2(\R/\Z)$, such that $u_1(x_2)$
coincides, near $x_2=0$, with the function $\sin\frac{1}{x_2}$, and
$u_3(x_1)$  coincides, near $x_1=0$, with the function
${\mathrm{sgn}} (x_1)$. Then by virtue of Theorem \ref{energyconth}
the shear flow
$$
u(x,t)=(u_1(x_2),0,u_3(x_1-tu_1(x_2))
$$
is a weak solution of the $3d$ Euler equations which conserves the
energy and which does not satisfy the hypothesis that are given in
\cite{shvydkoykh}.


\section{Conclusion}
We have used  the simplest example of a genuinely $3d$ flow to
obtain the following observations concerning the Euler equations:

\begin{enumerate}
\item In the class of H\"older spaces the space $C^1$ is the critical
space for the initial value problem of the $3d$ Euler equations  to
be locally, in time, well-posed {\it in the sense of Hadamard}. Old
and classical results \cite{lichtenstein} (see also
\cite{BealeKatoMajda} and \cite{Majda-Bertozzi}) have shown that the
$3d$ Euler equations are well-posed in $C^{1,\alpha}$, for every
$\alpha\in (0,1]$,  while  we have shown in section
\ref{loss-smoothness-section} that the $3d$ Euler equations are not
well-posed in $C^\beta$, for any $\beta\in (0,1)$. This observation
is also in agreement with the recent  result of Pak and Park
\cite{Pakpark}, who have established the local  well-posed, of the
$3d$ Euler equations, in the Besov space $B^{1}_{\infty,1}$. The
consistency between our result and that of \cite{Pakpark} is clear
from the   inclusion relations $C^{1,\alpha}\subset
B^{1}_{\infty,1}\subset C^1\,.$ Moreover, and as we have noted in
Remark \ref{Lemarie}, the analysis in section
\ref{loss-smoothness-section} can be adapted in more exotic spaces,
namely,  the shear  flow solutions, (\ref{shear}), of the $3d$ Euler
equations will provide examples of instabilities (i.e., the Cauchy
problem is not well-posed) in the in the Besov space $
B^{1}_{\infty,\infty}$, and in the Triebel-Lizorkin space
$F^1_{\infty,2}\,$.

\item The Kelvin-Helmholtz problem refers to  a free boundary problem
where in the $2d$ case limited regularity implies analyticity. We
show that in  $3d$, for closely related problems constructed with
the shear flow, this property is no more true. We propose an
explanation for this striking difference between the $2d$ and $3d$
case. This explanation is based on the fact that the linearized
operator of the Kelvin-Helmholtz problem is no longer
 elliptic in $3d$ as the situation  is in the
$2d$ case.

\item The relation between dissipation of energy and loss
of regularity is an essential issue in the statistical theory of
turbulence, in relation with the Kolmogorov Obukhov law. It has been
shown in the deterministic framework that a regularity of this type
implies conservation of energy. With the shear flow example we have
shown that there is no hope for a converse statement (even in the
case of solutions singular on a {\it slit} as in \cite{shvydkoykh}).
We observe that  in \cite{DeLellisandSzekelyhidi} De Lellis and
Szekelyhidi  constructed (see   cf. Theorem 1.1 a  ) an infinite set
of weak solutions
$$u \in C (\R_t; L^2(\R^3))$$ which satisfy both the strong and local energy equality (in the sense of
Definition 2.4 of   \cite{DeLellisandSzekelyhidi}) hence conserve energy.

The above observations may not invalidate the common {\it physical
belief} because the Kolmogorov Obukhov law belongs to the
statistical theory of turbulence, where statements and results are
true in some {\it averaged sense} . On the other hand, our family of
shear flow examples are genuinely laminar and therefore not
``turbulent." They are particular enough to be of measure zero with
respect to any reasonable {\it ensemble measure}    compatible with
the statistical theory of ideal (inviscid) turbulent flows (let us
recall that, to the best of our knowledge, no such measure has been
constructed, up to now, with full mathematical rigor).

Eventually the construction of \cite{DeLellisandSzekelyhidi}
involves limit of oscillating solutions and therefore is not
explicit but closer to the intuition of turbulence. It also relies
on the Baire category theorem. Hence it generates a residual set of
solutions which is dense in $C (\R_t; L_{\rm weak}^2(\R^3)).$  A
tentative justification of  the fact that `in the statistical theory
of turbulence conservation of energy may in general imply some
regularity of the underlined solutions' would be similar to the
situation in classical analysis theory where  a dense set may well
be a set of measure zero.

\end{enumerate}

\section*{Acknowledgments}
Claude Bardos would like to thank  the University of California,
Irvine and the Weizmann Institute of Science for their warm
hospitality where part of this work was completed.
The authors wish also to thank the anonymous referee for comments leading to the present improved version and the organizers of the Oberwolfach Workshop ID 0930 on Navier-Stokes equations (July 2009) where the material of this contribution was presented and discussed.
This work was
supported in part by the NSF grant no.~DMS-0708832, and by  the ISF
grant no. 120/06.

\end{document}